\documentclass[a4paper,12pt]{amsart}
\usepackage{amssymb}
\usepackage{ifthen}
 \usepackage[dvips]{graphicx}
\nonstopmode \numberwithin{equation}{section}
\usepackage{amssymb}
\usepackage{ifthen}
\usepackage{graphicx}
\usepackage{amsmath}
\usepackage[T1]{fontenc} 

\nonstopmode \numberwithin{equation}{section}
\theoremstyle{plain}
\newtheorem{thm}[equation]{Theorem}
\newtheorem{cor}[equation]{Corollary}
\newtheorem{lem}[equation]{Lemma}
\newtheorem{prop}{Proposition}

\newtheorem{conj}{Conjecture}

\theoremstyle{definition}
\newtheorem{defn}{Definition}[section]

\newtheorem{prob}{Problem}
\newtheorem{rem}{Remark}[section]

\newcounter{minutes}
\divide\time by 60
\newcounter{hours}
\multiply\time by 60
\addtocounter{minutes}{-\time}

\newcounter {own}
\def\theown {\thesection       .\arabic{own}}

\newenvironment{pf}[1][]{%
 \vskip 3mm
 \noindent
 \ifthenelse{\equal{#1}{}}%
  {{\slshape Proof. }}%
  {{\slshape #1.} }%
 }%
{\qed\bigskip}

\newcounter{alphabet}
\newcounter{tmp}

\newcommand{\real}{{\operatorname{Re}\,}}

\def\be{\begin{equation}}
\def\ee{\end{equation}}

\newcommand{\bee}{\begin{enumerate}}
\newcommand{\eee}{\end{enumerate}}

\newcommand{\blem}{\begin{lem}}
\newcommand{\elem}{\end{lem}}
\newcommand{\bthm}{\begin{thm}}
\newcommand{\ethm}{\end{thm}}
\newcommand{\bcor}{\begin{cor}}
\newcommand{\ecor}{\end{cor}}
\newcommand{\beg}{\begin{examp}}
\newcommand{\eeg}{\end{examp}}
\newcommand{\begs}{\begin{examples}}
\newcommand{\eegs}{\end{examples}}

\newcommand{\bdefn}{\begin{defn}}
\newcommand{\edefn}{\end{defn}}

\newcommand{\bprob}{\begin{prob}}
\newcommand{\eprob}{\end{prob}}
\newcommand{\bei}{\begin{itemize}}
\newcommand{\eei}{\end{itemize}}

\newcommand{\bcon}{\begin{conj}}
\newcommand{\econ}{\end{conj}}
\newcommand{\bcons}{\begin{conjs}}
\newcommand{\econs}{\end{conjs}}
\newcommand{\bprop}{\begin{prop}}
\newcommand{\eprop}{\end{prop}}
\newcommand{\br}{\begin{rem}}
\newcommand{\er}{\end{rem}}
\newcommand{\brs}{\begin{rems}}
\newcommand{\ers}{\end{rems}}
\newcommand{\bo}{\begin{obser}}
\newcommand{\eo}{\end{obser}}
\newcommand{\bos}{\begin{obsers}}
\newcommand{\eos}{\end{obsers}}
\newcommand{\bpf}{\begin{pf}}
\newcommand{\epf}{\end{pf}}
\newcommand{\ba}{\begin{array}}
\newcommand{\ea}{\end{array}}
\newcommand{\beq}{\begin{eqnarray}}
\newcommand{\beqq}{\begin{eqnarray*}}
\newcommand{\eeq}{\end{eqnarray}}
\newcommand{\eeqq}{\end{eqnarray*}}

\begin{document}

\title{Bohr phenomenon for certain   Subclasses of Harmonic  Mappings}

\author{Vasudevarao Allu}
\address{Vasudevarao Allu,
School of Basic Science,
Indian Institute of Technology Bhubaneswar,
Bhubaneswar-752050, Odisha, India.}
\email{avrao@iitbbs.ac.in}

\author{Himadri Halder}
\address{Himadri Halder,
School of Basic Science,
Indian Institute of Technology Bhubaneswar,
Bhubaneswar-752050, Odisha, India.}
\email{hh11@iitbbs.ac.in}

\subjclass[{AMS} Subject Classification:]{Primary 30C45, 30C50, 30C80}

\def\thefootnote{}
\footnotetext{ {\tiny File:~\jobname.tex,
printed: \number\year-\number\month-\number\day,
          \thehours.\ifnum\theminutes<10{0}\fi\theminutes }
} \makeatletter\def\thefootnote{\@arabic\c@footnote}\makeatother

\begin{abstract}
 The Bohr phenomenon for analytic functions of the form $f(z)=\sum_{n=0}^{\infty} a_{n}z^{n}$, first introduced by Harald Bohr in 1914, deals with finding the largest radius 
 $r_{f}$, $0<r_{f}<1$, such that the inequality $\sum_{n=0}^{\infty} |a_{n}z^{n}| \leq 1$ holds whenever the inequality $|f(z)|\leq 1 $ holds in the unit disk 
 $\mathbb{D}=\{z \in \mathbb{C}: |z|<1\}$. The exact value of this largest radius known as Bohr radius, which has been established to be  $r_{f}=1/3$. 
The Bohr phenomenon \cite{Abu-2010} for harmonic functions  $f$ of the form  $f(z)=h(z)+\overline {g(z)}$, where $h(z)=\sum_{n=0}^{\infty} a_{n}z^{n}$ 
and $g(z)=\sum_{n=1}^{\infty} b_{n}z^{n}$ is to find the largest radius $r_{f}$, $0<r_{f}<1$ such that 
$$\sum\limits_{n=1}^{\infty} (|a_{n}|+|b_{n}|) |z|^{n}\leq d(f(0),\partial f(\mathbb{D})) 
$$ 
holds for $|z|\leq r_{f}$, here $d(f(0),\partial f(\mathbb{D})) $ denotes the Euclidean distance between $f(0)$ and the boundary of $f(\mathbb{D})$. 
In this paper, we investigate the Bohr radius for several classes of harmonic functions in the unit disk $\mathbb{D}.$
\end{abstract}

\maketitle
\pagestyle{myheadings}
\markboth{Vasudevarao Allu and  Himadri Halder}{Bohr phenomenon for certain   Subclasses of Harmonic Mappings}

\section{Introduction and Preliminaries}

Let $\mathcal{F}$ be the class of  analytic functions $f$ on the  unit disk $\mathbb{D}=\{z \in \mathbb{C}: |z|<1\}$ with $|f(z)|\leq 1 $ in $\mathbb{D}$. 
Each  function  $f \in \mathcal{F}$ has the following series representation  
\begin{equation} \label{himadri-vasu-p1-e-1.1}
f(z)=\sum\limits_{n=0}^{\infty} a_{n}z^{n}
\end{equation}
and the power series \eqref{himadri-vasu-p1-e-1.1} converges uniformly on every compact subset of $\mathbb{D}$.  For a  given power series of
the form \eqref{himadri-vasu-p1-e-1.1}, its majorant series  
is defined by
\begin{equation} \label{himadri-vasu-p1-e-1.2}
M_{f}(r)=\sum\limits_{n=0}^{\infty} |a_{n}z^{n}|= \sum\limits_{n=0}^{\infty} |a_{n}|r^{n}
\end{equation}
for $|z|=r <1.$
By the series comparison test, it is easy to see that the series \eqref{himadri-vasu-p1-e-1.1} and \eqref{himadri-vasu-p1-e-1.2} converge or diverge together in $\mathbb{D}$.  
Clearly, for each $f\in \mathcal{F}$ we denote $M(r):= M_{f}(r)$ which
is an increasing function for $0\leq r<1$ and  $M_{f}(0)=|a_{0}|=|f(0)| \leq 1.$ 
We note that for certain functions $f$ in $\mathcal{F}$, the majorant series $M_{f}(r)$ takes the values larger than $1$ (see \cite{h p bose}). 
Therefore, it is natural to ask for the possible values of $r\in [0,1)$ so that 

\begin{equation} \label{himadri-vasu-p1-e-1.3}
M_{f}(r)=\sum\limits_{n=0}^{\infty} |a_{n}|r^{n} \leq 1
\end{equation}
holds for all $f\in \mathcal{F}.$  In 1914, Harald Bohr \cite{Bohr-1914} observed that the inequality  \eqref{himadri-vasu-p1-e-1.3} is true for  $|z| \leq 1/6 $ and this practice 
was further expanded by Wiener, Riesz, and Schur who independently established the inequality \eqref{himadri-vasu-p1-e-1.3} on the disk $|z|=r \leq 1/3$ and that $1/3$, now 
called Bohr radius, is the best possible  (see \cite{paulsen-2002,sidon-1927,tomic-1962}). The inequality \eqref{himadri-vasu-p1-e-1.3} on $[0,1)$ is usually known 
as Bohr inequality for the class $\mathcal{F}.$  

\vspace{8mm}
\hspace{5mm} We note that the inequality \eqref{himadri-vasu-p1-e-1.3}, in turn,  equivalent to
\begin{equation} \label{himadri-vasu-p1-e-1.4}
\sum\limits_{n=1}^{\infty} |a_{n}z^{n}| \leq 1-|a_{0}|=d(f(0),\partial f(\mathbb{D}) ) \quad \mbox{for} \quad |z|\leq \dfrac{1}{3},
\end{equation}
where $d(f(0),\partial f(\mathbb{D}) )$ denotes the Euclidean distance between $f(0)$ and the boundary of $f(\mathbb{D})$. Functions $f \in \mathcal{F}$ 
satisfying  \eqref{himadri-vasu-p1-e-1.4},  sometimes are said to satisfy the classical Bohr's phenomenon. It is important to note that the existance of the 
radius $1/3$ in \eqref{himadri-vasu-p1-e-1.4} is independent of the coefficients of the power series \eqref{himadri-vasu-p1-e-1.1} {\it i.e.,} in a better way we 
can demonstrate this fact by saying that a Bohr phenomenon appears in the class of analytic self-maps of the unit disk $\mathbb{D}$. Later, the concept of 
Bohr radius has been developed to the class of analytic functions from $\mathbb{D}$ into certain domains $G \subseteq \mathbb{C}$ and it is found that the 
above radius varies with respect to different types of domains in $\mathbb{C}$. Therefore, it is natural to generalize the above Bohr phenomenon as: for a given 
domain $G \subseteq \mathbb{C}$, to find the largest radius $r_{G}>0$ such that 
\begin{equation} \label{himadri-vasu-p1-e-1.5}
\sum\limits_{n=1}^{\infty} |a_{n}z^{n}| \leq d(f(0),\partial f(\mathbb{D}) ) \quad \mbox{in} \quad |z|<r_{G}
\end{equation}
for all functions $f$ belong  to the class of analytic functions in $\mathbb{D}$ such that $f(\mathbb{D}) \subseteq G$. 
For a convex domain $G \subseteq \mathbb{D}$, Aizenberg \cite{aizn-2007} has obtained the largest radius $r_{G}$, for which the 
inequality  \eqref{himadri-vasu-p1-e-1.5} holds, coincides with the classical Bohr radius $1/3$ which cannot be improved further; while 
Abu-Muhanna \cite{Abu-2010} has proved that the same inequality holds for the largest radius $r_{G}=3-2\sqrt{2}$ for any proper simply connected domain $G$.

\vspace{5 mm}

In \cite{aizenberg-2001}, the existance of Bohr phenomenon for the class of holomorphic functions on complex manifold has been established. 
In \cite{abu-2014}, the Bohr phenomenon for mappings from $\mathbb{D}$ into the wedge-domain 
$W_{\alpha}=\{w: |arg \, w|<\pi \alpha /2, 1\leq \alpha \leq 2\}$ has been studied. A simple observation shows that $W_{\alpha}$ is convex only for $\alpha = 1$, which 
in fact coincides with the right half-plane. In 2018, Bhowmik and Das \cite{bhowmik-2018} studied the Bohr radius for the class of concave univalent functions 
with openinig angle $\pi \alpha$, $1\leq \alpha \leq 2$ and also obtained the Bohr radius for the class of starlike functions of order 
$\alpha$, for $0\leq \alpha \leq 1/2$. In \cite{Ali-2019}, Ali {\it et al.} obtained the Bohr radius for the class of convex functions of order $\alpha,$ for $-1/2<\alpha <1$.
 Improved Bohr radius for the classes of starlike and convex functions of order $\alpha$ with negative coefficients have also been 
 obtained in \cite{Ali-2019}. In 2017, Ali {\it et al.} \cite{r m ali jmma} obtained the Bohr radius for $n$-$symmetric$ and alternating series as well as for the class of 
even analytic functions. The Bohr radius for various classes of functions {\it e.g.} locally univalent harmonic mappings, $k$-quasiconformal mappings, bounded 
harmonic functions, lacunary series have been studied extensively  in \cite{ponnusamy-locally univalent,ponnusamy-jmma}. 
For more intriguing aspects of Bohr phenomenon we refer to \cite{aizn-2000,seraj-2019,bene-2004,evdoridis-2019,kayumova-2018,Z-Liu-2019} and the references therein. \\

We now present the definition of subordination for analytic functions which is useful to prove some of our theorems. Let $f$ and $g$ be two analytic functions in $\mathbb{D}$.  
Then $f$ is subordinate to $g$, written $f \prec g$, if there exists an anlytic function $\omega$ in $\mathbb{D}$ with $\omega(0)=0$ and $|\omega(z)|<1$, such that $f(z)=g(\omega(z)).$
In particular, if $g$ is univalent in $\mathbb{D},$ then $f$ is subordinate to $g$ provided $f(0)=g(0)$ and $f(\mathbb{D}) \subseteq g(\mathbb{D}).$\\

Let $\mathcal{H}$ be the class of all complex-valued harmonic functions $f=h+\overline{g}$ defined on $\mathbb{D}$, where $h$ and $g$ are analytic 
in $\mathbb{D}$ with the normalization $h(0)=h'(0)-1=0$ and $g(0)=0$. Let 
$$
\mathcal{H}_{0}=\{f=h+\overline{g}\in \mathcal{H}:g'(0)=0\}.
$$ 
Then, each $f=h+\overline{g}\in \mathcal{H}_{0}$ has the following form 
\begin{equation} \label{himadri-vasu-p1-e-1.6}
h(z)=z+ \sum \limits _{n=2}^{\infty} a_{n}z^{n} \quad \mbox{and} \quad g(z)=\sum\limits_{n=2}^{\infty} b_{n}z^{n}.
\end{equation}
In 2013, Ponnusamy {\it et al.} \cite{S.Ponnusamy variability regions-2013}  considered 
the following classes   
$$\mathcal{P}_{\mathcal{H}}=\{f=h+\bar{g} \in \mathcal{H} : \real h'(z)>|g'(z)| \quad  \mbox{for} \quad z \in \mathbb{D}\}
$$
and $\mathcal{P}^{0}_{\mathcal{H}}=\mathcal{P}_{\mathcal{H}} \cap \mathcal{H}_{0}.$ It is known that  functions in  $\mathcal{P}_{\mathcal{H}}$ 
are close-to-convex (see \cite{S.Ponnusamy variability regions-2013}). Motivated by the above classes, Li and Ponnusamy \cite{Injectivity section} have studied the following classes 
$$\mathcal{P}_{\mathcal{H}}(\alpha)=\{f=h+\overline{g} \in \mathcal{H} : \real (h'(z)-\alpha)>|g'(z)| \,  \quad \mbox{for} \quad  z \in \mathbb{D}\}
$$
and $\mathcal{P}^{0}_{\mathcal{H}}(\alpha)=\mathcal{P}_{\mathcal{H}}(\alpha) \cap \mathcal{H}_{0}.$  Clearly, 
$\mathcal{P}_{\mathcal{H}}(\alpha) \subseteq \mathcal{P}_{\mathcal{H}}$ and $\mathcal{P}^{0}_{\mathcal{H}}(\alpha) \subseteq \mathcal{P}^{0}_{\mathcal{H}}$ for $0\leq \alpha <1$. \\

In \cite{Injectivity section}, Li and Ponnusamy have 
proved that  functions in $\mathcal{P}^{0}_{\mathcal{H}}(\alpha)$ are univalent for $0\leq \alpha <1$. Further, the coefficient bounds and the univalency of sections for functions in 
$\mathcal{P}^{0}_{\mathcal{H}}(\alpha)$ have been studied in \cite{Injectivity section}.

\vspace{3mm}
In 2016, Li and Ponnusamy \cite{Li and Ponnysamy -2016} studied the following class 
$$
\widetilde{\mathcal{G}}^{0}_{\mathcal{H}}(\beta)=\left\{f=h+\overline{g} \in \mathcal{H}_{0}: \real \left(\dfrac{h(z)}{z} \right)-\beta > \left|\dfrac{g(z)}{z} \right| \, \mbox{for} \, z\in \mathbb{D}\right\}.
$$
It is known that $\widetilde{\mathcal{G}}^{0}_{\mathcal{H}}(\beta) \subseteq \widetilde{\mathcal{G}}^{0}_{\mathcal{H}}(0)$ for $0\leq \beta <1.$ It is proved that the 
harmonic convolution of functions in  $\mathcal{P}^{0}_{\mathcal{H}}(\beta)$ with functions in $\widetilde{\mathcal{G}}^{0}_{\mathcal{H}}(\beta)$ is 
univalent and close-to-convex harmonic in $\mathbb{D}$ with  certain conditions on the parameters $\alpha$ and $\beta$. 

\vspace{3mm}

In 1977, Chichra \cite{chichra-1977} first introduced the class $\mathcal{W}(\alpha)$, consisting of normalized analytic functions $h$, satisfying the condition 
$\real (h'(z)+\alpha z h''(z))>0$ for $ z\in \mathbb{D}$ and $\alpha \geq 0.$ Functions in $\mathcal{W}(\alpha)$ are univalent. 
Further, Chichra \cite{chichra-1977} has shown that   functions in the  class $\mathcal{W}(\alpha)$ constitute a subclass of close-to-convex functions in $\mathbb{D}$. 
In 2014, Nagpal and Ravichandran \cite{nagpal-2014} studied the following class   
$$\mathcal{W}^{0}_{\mathcal{H}}=\{f=h+\overline{g} \in \mathcal{H}: \quad \real (h'(z)+ z h''(z))>|g'(z)+zg''(z)| \quad \mbox{for} \quad z\in \mathbb{D}\}
$$  
and obtained the coefficient bounds for functions in  $\mathcal{W}^{0}_{\mathcal{H}}$.
Motivated by the class $\mathcal{W}^{0}_{\mathcal{H}}$, Ghosh and Vasudevarao \cite{nirupam-2019} have  recently studied the following class 
$\mathcal{W}^{0}_{\mathcal{H}}(\alpha),$ where 
$$\mathcal{W}^{0}_{\mathcal{H}}(\alpha)=\{f=h+\overline{g} \in \mathcal{H}_{0}:\,  \real (h'(z)+\alpha z h''(z))>|g'(z)+\alpha zg''(z)| \quad \mbox{for} \quad  z\in \mathbb{D} \}.
$$  

\vspace{3mm}

In 1977, Chichra \cite{chichra-1977} studied the class $\mathcal{G}(\alpha)$ consists of analytic functions $f\in \mathcal{A}$ satisfying the condition 
$$
\real \left((1-\alpha) \dfrac{f(z)}{z}+\alpha f'(z)\right)>0
$$ 
for  $z \in \mathbb{D}$ and $\alpha \geq 0,$ where  $\mathcal{A}$ is the class of analytic functions $h$ with the normalization $h(0)=h'(0)-1=0.$ Chichra \cite{chichra-1977} 
has proved that functions in $\mathcal{G}(\alpha)$ are univalent in 
$\mathbb{D}$ for $\alpha \geq 1$, while  functions in $\mathcal{G}(0)$ are univalent in $|z|<\sqrt{2}-1$. Liu \cite{Liu-2000} has shown that, if 
$0\leq \beta < \alpha $ then  $\mathcal{G}(\alpha) \subseteq \mathcal{G}(\beta)$.  
Motivated by the above class, Liu {\it et al.} \cite{Li-Mei-2018} have introduced the following  subclasses of $\mathcal{H}_{0}$  
\begin{eqnarray*}
	\mathcal{G}_{\mathcal{H}}^{k}(\alpha):&=& \left\{ f=h+\overline{g} \in \mathcal{H}_{0} ^{k}  :  \, 
	\real \left((1-\alpha)\dfrac{h(z)}{z} + \alpha h'(z)\right)\right.\\[2mm]
	&& \hspace*{4cm} > \left. \left|(1-\alpha)\dfrac{g(z)}{z} + \alpha g'(z) \right| \quad\mbox{ for }  z\in \mathbb{D}\right\} 
\end{eqnarray*}
where 
$\mathcal{H}_{0} ^{k}=\{f=h+\overline{g} \in \mathcal{H}: h'(0)-1=g'(0)=h''(0)=\cdots =h^{k}(0)=g^(k)(0)=0\}$ 
and $\mathcal{H}_{0}^{1}\equiv \mathcal{H}_{0}$ for some $\alpha \geq 0,$ $ k\geq 1$. Thus, every $f=h+\overline{g} \in \mathcal{H}_{0} ^{k}$ has the following representation 
\begin{equation} \label{himadri-vasu-p1-e-1.7}
h(z)=z+\sum\limits_{n=k+1}^{\infty} a_{n} z^{n} \quad \mbox{and} \quad g(z)=\sum\limits_{n=k+1}^{\infty} b_{n} z^{n}.
\end{equation}
Note that 
$\mathcal{G}_{\mathcal{H}}^{n}(\alpha) \subseteq \mathcal{G}_{\mathcal{H}}^{k}(\alpha) \subseteq \mathcal{G}_{\mathcal{H}}^{1}(\alpha)$
for each 
$n \geq k \geq 1$ and $\mathcal{G}_{\mathcal{H}}^{1}(1)=\mathcal{P} _{\mathcal{H}}.$ In \cite{Li-Mei-2018}, it has been proved that functions in 
$\mathcal{G}_{\mathcal{H}}^{1}(\alpha)$ are univalent in $\mathbb{D}$ and functions in $\mathcal{G}_{\mathcal{H}}^{1}(0)$ are starlike in the disk $|z|<\sqrt{2}-1$. 
The sections of functions  and many other geometric properties of functions in $\mathcal{G}_{\mathcal{H}}^{k}(\alpha)$ have also been established in \cite{Li-Mei-2018}.

\vspace{3mm}
The following two subclasses of harmonic functions have been introduced by Ghosh and Vasudevarao \cite{nirupam cvee,nirupam-bull aus math},
$$\mathcal{B}^{0}_{\mathcal{H}}(M)=\{f=h+\overline{g} \in \mathcal{H}_{0}: |zh''(z)|\leq M-|zg''(z)|, \quad z \in \mathbb{D}\quad \mbox{and } \quad M>0\},
$$
$$\mathcal{P}^{0}_{\mathcal{H}}(M)=\{f=h+\overline{g} \in \mathcal{H}_{0}: \real (zh''(z))> -M+|zg''(z)|, \quad z \in \mathbb{D} \quad \mbox{and } \quad M>0\}.
$$ 
The analytic and geometric properties of the classes $\mathcal{B}^{0}_{\mathcal{H}}(M)$ and $\mathcal{P}^{0}_{\mathcal{H}}(M)$ have been 
extensively studied in \cite{nirupam cvee,nirupam-bull aus math}. The subclasses 
$\mathcal{B}^{0}_{\mathcal{H}}(M)$ and $\mathcal{P}^{0}_{\mathcal{H}}(M)$ are not only the generalizations of analytic functions 
but also they are closely related to the analytic subclasses $\mathcal{B}(M)$ and $\mathcal{P}(M)$ respectively, where the analytic subclasses are defined by 
$$\mathcal{B}(M)=\{\phi \in \mathcal{A}:|z\phi ''(z)|\leq M , \quad z\in \mathbb{D} \quad \mbox{and } \quad M>0\},
$$ 
$$\mathcal{P}(M)=\{\phi \in \mathcal{A}: \real (z\phi ''(z))> -M, \quad  z\in \mathbb{D} \quad \mbox{and} \quad M>0\}.
$$
The classes $\mathcal{B}(M)$ and $\mathcal{P}(M)$ have been studied extensively by Mocanu \cite{Mocanu-1992} and Ponnusamy \cite{ponnusamy -2002}. In 2002, Ponnusamy {\it et al.} \cite{ponnusamy -2002} proved that functions in  $\mathcal{B}(M)$ are univalent and starlike for $0<M\leq 1$, and convex for $0<M\leq 1/2$ where the authors in \cite{Ali 1995} have shown that each function in $\mathcal{P}(M)$ 
is univalent and starlike for $0<M<1/(ln \, 4).$ The regions of variability for the classes $\mathcal{B}(M)$ and $\mathcal{P}(M)$ have been studied by Ponnusamy {\it et al.} \cite{S.Ponnusamy variability regions-2009}.

\vspace{4mm}

We define the class $\mathcal{T}_{\mathcal{B}^{0}_{\mathcal{H}}(M)}$ denote the class of functions $f \in \mathcal{H}_{0}$ of the form \eqref{himadri-vasu-p1-e-1.6} satifying 
$\sum_{n=2}^{\infty} n(n-1) (|a_{n}|+|b_{n}|) \leq M, $ for some $M>0.$ In view of \cite[Theorem 2.4]{nirupam cvee}, we have 
$\mathcal{T}_{\mathcal{B}^{0}_{\mathcal{H}}(M)} \subseteq \mathcal{B}^{0}_{\mathcal{H}}(M).$  

\begin{lem} \label{himadri-vasu-p1-lemma-001}  \cite{Injectivity section}
Let $f \in \mathcal{P}^{0}_{\mathcal{H}}(\alpha) $ be given by \eqref{himadri-vasu-p1-e-1.6}. Then for any $n \geq 2$, 
\begin{enumerate}
	\item[(i)] $\displaystyle |a_n| + |b_n|\leq \frac {2(1-\alpha)}{n}; $\\[2mm]
	
	\item[(ii)] $\displaystyle ||a_n| - |b_n||\leq \frac {2(1-\alpha)}{n};$\\[2mm]
	
	\item[(iii)] $\displaystyle |a_n|\leq \frac {2(1-\alpha)}{n}.$
\end{enumerate}
All the results are sharp, with $f(z)=(1-\alpha)(-z-2\, \mbox{log}(1-z))+\alpha z$ being the extremal.
\end{lem}

The following result shows that functions in $\mathcal{W}^0_{\mathcal{H}}(\alpha)$ are univalent for $\alpha \geq 0,$ and they are closely related to functions in $\mathcal{W}(\alpha).$
\begin{lem} \label{himadri-vasu-p1-lemma-002} \cite{nirupam-2019}
	The harmonic mapping $f = h + \overline{g}$ belongs to $\mathcal{W}^0_{\mathcal{H}}(\alpha)$ if, and only if, the analytic function $F = h + \epsilon {g}$ belongs to $\mathcal{W}(\alpha)$ for each $|\epsilon| = 1$.
\end{lem}
The coefficient bounds and the sharp growth estimates for  functions in the class  $\mathcal{W}^{0}_{\mathcal{H}}(\alpha)$ have been studied in \cite{nirupam-2019}.

\begin{lem} \label{himadri-vasu-p1-lemma-003} \cite{nirupam-2019}
	Let $f \in \mathcal{W}^{0}_{\mathcal{H}}(\alpha)$ for $\alpha \geq 0$ and be of the form \eqref{himadri-vasu-p1-e-1.6}. 
	Then for $n\geq 2,$
	\begin{enumerate}
		\item[(i)] $\displaystyle |a_n| + |b_n|\leq \frac {2}{\alpha n^2 + n(1 - \alpha)}; $\\[2mm]
		
		\item[(ii)] $\displaystyle ||a_n| - |b_n||\leq \frac {2}{\alpha n^2 + n(1 - \alpha)};$\\[2mm]
		
		\item[(iii)] $\displaystyle |a_n|\leq \frac {2}{\alpha n^2 + n(1 - \alpha)}.$
	\end{enumerate}
	All these results are sharp for the function $f$ given by $f(z) = z +  \sum _{n = 2}^{\infty}\dfrac{2}{\alpha n^2 + n(1 - \alpha)}z^n.$
\end{lem}

\begin{lem} \label{himadri-vasu-p1-lemma-003-a}\cite{nirupam-2019}
Let $f = h + \overline{g} \in \mathcal{W}^0_{\mathcal{H}}(\alpha)$  be given by \eqref{himadri-vasu-p1-e-1.6} with  $0<\alpha \leq 1.$ Then 
\begin{equation} \label{himadri-vasu-p1-e-1.12}
|z| + 2\sum_{n = 2}^{\infty}\frac{(-1)^{n-1} |z|^n}{\alpha n^2 + n(1 - \alpha)} \leq |f(z)|\leq |z| + 2\sum_{n = 2}^{\infty}\frac{|z|^n}{\alpha n^2 + n(1 - \alpha)}.
\end{equation}
Both inequalities are  sharp when $f$ is given by $f(z) = z + \sum _{n = 2}^{\infty}\frac{2}{\alpha n^2 + n(1 - \alpha)}z^n$, or  its rotations.
\end{lem}
The following lemmas are usefull to establish the Bohr radius for the class $\mathcal{G}_{\mathcal{H}}^{k}(\alpha).$ 

\begin{lem} \label{himadri-vasu-p1-lemma-004} \cite{Li-Mei-2018}
	Let $f=h+\overline{g} \in \mathcal{G}_{\mathcal{H}}^{k}(\alpha) $ where $h$ and $g$ be given by \eqref{himadri-vasu-p1-e-1.7} with $k\geq 1$. Then for any $n\geq k+1,$ 
	\begin{enumerate}
		\item[(i)] $\displaystyle |a_n| + |b_n|\leq \frac {2}{1 + (n-1)\alpha}; $\\[2mm]
		
		\item[(ii)] $\displaystyle ||a_n| - |b_n||\leq \frac {2}{1 + (n-1)\alpha};$\\[2mm]
		
		\item[(iii)] $\displaystyle |a_n|\leq \frac {2}{1 + (n-1)\alpha}.$
	\end{enumerate} 
	All the results are sharp, with the functions $f_{i}(z)=z+ \sum_{j=1}^{\infty} \dfrac{2}{1+ji \alpha } z^{ji+1}, \quad i=k,k+1,...,2k-1$ being the extremals. 
\end{lem} 
\begin{lem} \label{himadri-vasu-p1-lemma-005} \cite{Li-Mei-2018}
	Let $f=h+\bar{g} \in \mathcal{G}_{\mathcal{H}}^{k}(\alpha) $ as in the form \eqref{himadri-vasu-p1-e-1.7} with $\alpha \geq 0.$ Then 
	\begin{equation} \label{himadri-vasu-p1-e-1.16}
	|z|+2\sum\limits_{n=1}^{\infty} \dfrac{(-1)^{n}}{1+nk \alpha } |z|^{nk+1} \leq |f(z)| \leq |z|+2\sum\limits_{n=1}^{\infty} \dfrac{1}{1+nk \alpha } |z|^{nk+1}.
	\end{equation}
	This result is sharp with the function $f(z)=z+ \sum_{n=1}^{\infty} \dfrac{2}{1+nk \alpha } z^{nk+1}$ and its rotations being the extremals.
\end{lem}
The following one-to-one correspondence between the classes $\mathcal{B}^{0}_{\mathcal{H}}(M)$($\mathcal{P}^{0}_{\mathcal{H}}(M)$ respectively) and 
$\mathcal{B}(M)$($\mathcal{P}(M)$ respectively) obtained by Ghosh and Vasudevarao \cite{nirupam cvee,nirupam-bull aus math}. 
\begin{lem} \label{himadri-vasu-p1-lemma-006} \cite{nirupam cvee}
	A harmonic function $f=h+\overline{g} $ is in $\mathcal{B}^{0}_{\mathcal{H}}(M)$ if, and only if, $F_{\epsilon}=h+\epsilon g$ is in $\mathcal{B}(M)$ for each $\epsilon (|\epsilon|=1).$
\end{lem}
\begin{lem} \label{himadri-vasu-p1-lemma-007} \cite{nirupam-bull aus math}
A harmonic function $f=h+\overline{g} $ is in $\mathcal{P}^{0}_{\mathcal{H}}(M)$ if, and only if, $F_{\epsilon}=h+\epsilon g$ is in $\mathcal{P}(M)$ for each $\epsilon (|\epsilon|=1).$
\end{lem}
\begin{lem} \cite{nirupam cvee} \label{himadri-vasu-p1-lemma-008}
	Let $f=h+\bar{g}\in \mathcal{B}^{0}_{\mathcal{H}}(M)$ for some $M>0.$ Then 
	\begin{equation} \label{himadri-vasu-p1-e-1.20}
	|z|-\dfrac{M}{2}|z|^{2} \leq |f(z)| \leq |z|+\dfrac{M}{2}|z|^{2}.
	\end{equation}
	Both inequalities are sharp.
\end{lem}

\begin{lem} \label{himadri-vasu-p1-lemma-009} \cite{nirupam-bull aus math}
Let $f=h+\overline{g}\in \mathcal{P}^{0}_{\mathcal{H}}(M)$ for some $M>0$ be of the form \eqref{himadri-vasu-p1-e-1.6}. Then for $n\geq 2,$ 
\begin{enumerate}
	\item[(i)] $\displaystyle |a_n| + |b_n|\leq \frac {2M}{n(n-1)}; $\\[2mm]
	
	\item[(ii)] $\displaystyle ||a_n| - |b_n||\leq \frac {2M}{n(n-1)};$\\[2mm]
	
	\item[(iii)] $\displaystyle |a_n|\leq \frac {2M}{n(n-1)}.$
\end{enumerate}
The results are sharp for the function $f$ given by 
$f'(z)=1-2M\, ln\, (1-z) .$	
\end{lem}
In \cite{nirupam-bull aus math}, Nirupam and Vasudevarao have proved the folowing growth theorem  for the class $\mathcal{P}^{0}_{\mathcal{H}}(M)$  where  the  right hand inequality is sharp.
\begin{thm} \label{himadri-vasu-p1-theorem-006-b} \cite{nirupam-bull aus math}
	Let $f \in \mathcal{P}^{0}_{\mathcal{H}}(M)$. Then 
	\begin{equation} \label{himadri-vasu-p1-e-2.10-b}
	|z| -2M \sum\limits_{n=2}^{\infty} \dfrac{|z|^{n}}{n(n-1)} \leq |f(z)| \leq |z| + 2M \sum\limits_{n=2}^{\infty} \dfrac{|z|^{n}}{n(n-1)}.
	\end{equation}
	The right hand inequality is sharp for the function $f$ given by $f'(z)=1-2M\, \mbox{ln}\, (1-z).$
\end{thm}

\section{Main results}

 First we prove the following sharp growth estimate for the class $\mathcal{P}^{0}_{\mathcal{H}}(\alpha)$.

\begin{thm} (Growth estimate)\label{himadri-vasu-p1-theorem-001-a}
	Let $f=h+\overline{g} \in \mathcal{P}^{0}_{\mathcal{H}}(\alpha)$ with $0\leq \alpha <1$. Then 
	\begin{equation} \label{himadri-vasu-p1-e-2.2}
	|z|+ \sum\limits_{n=2}^{\infty}  \dfrac{2(1-\alpha)(-1)^{n-1}}{n} |z|^{n} \leq |f(z)| \leq |z|+ \sum\limits_{n=2}^{\infty}  \dfrac{2(1-\alpha)}{n} |z|^{n}.
	\end{equation} Both the inequalities are sharp.
\end{thm}

Using Lemma \ref{himadri-vasu-p1-lemma-001} and Theorem \ref{himadri-vasu-p1-theorem-001-a}, we obtain the Bohr radius for the class $\mathcal{P}^{0}_{\mathcal{H}}(\alpha).$ 

\begin{thm} \label{himadri-vasu-p1-theorem-001}
	Let $f \in \mathcal{P}^{0}_{\mathcal{H}}(\alpha)$ be given by \eqref{himadri-vasu-p1-e-1.6} with $0\leq \alpha <1$. 
	Then $$|z|+\sum\limits_{n=2}^{\infty} (|a_{n}|+|b_{n}|)|z|^{n} \leq d(f(0), \partial f (\mathbb{D}))$$ for $|z|=r\leq r_{f}$, where $r_{f}$ is the unique positive root of  
	\begin{equation} \label{himadri-vasu-p1-e-2.4}
	r+2(1-\alpha)\sum \limits _{n=2}^{\infty}\dfrac{r^{n}}{n}=1+2(1-\alpha)\sum \limits _{n=2}^{\infty}\dfrac{(-1)^{n-1}}{n} 
	\end{equation} 
	in $(0,1)$.	The radius $r_{f}$ is the best possible.
\end{thm}
\begin{rem}
If  $g \equiv 0$ and $\alpha=0$ then the class $\mathcal{P}^{0}_{\mathcal{H}}(\alpha)$ reduces to the class $\mathcal{P}$ of normalized analytic functions $h$ such that $\real h'(z)>0$ in $\mathbb{D}$ 
(see \cite{macgregor-1962}).
By taking  $\alpha=0$ in \eqref{himadri-vasu-p1-e-2.4}, we obtain the Bohr radius $r_{f} \approx 0.285194$ for the class $\mathcal{P}.$
\end{rem}
We prove the sharp growth estimate for functions in $\widetilde{\mathcal{G}}^{0}_{\mathcal{H}}(\beta)$ and using this we obtain 
 the Bohr phenomenon for the class $\widetilde{\mathcal{G}}^{0}_{\mathcal{H}}(\beta)$.  
\begin{thm} \label{himadri-vasu-p1-theorem-002}
	Let $f =h+\overline{g} \in \widetilde{\mathcal{G}}^{0}_{\mathcal{H}}(\beta)$ for $0\leq \beta < 1/2$ be given by \eqref{himadri-vasu-p1-e-1.6}. 
	Then 
	$$|z|+\sum\limits_{n=2}^{\infty}(|a_{n}|+|b_{n}|)|z|^{n} \leq d(f(0),\partial f(\mathbb{D}))
	$$ 
	for $|z|=r\leq r_{f}:=(-1-\beta +\sqrt{1+6\beta -7\beta ^{2}})/2(1-2\beta),$
	 where $r_{f}$ is the unique positive root of $$
	(1-2\beta)r^{2}+(1+\beta)r-\beta =0. $$ The radius $r_{f}$ is the best possible.
\end{thm}
Using Lemmas \ref{himadri-vasu-p1-lemma-003} and \ref{himadri-vasu-p1-lemma-003-a}, we find the sharp Bohr radius for the class $\mathcal{W}^{0}_{\mathcal{H}}(\alpha)$.
\begin{thm} \label{himadri-vasu-p1-theorem-003}
	Let $f \in \mathcal{W}^{0}_{\mathcal{H}}(\alpha)$  be given by \eqref{himadri-vasu-p1-e-1.6}. Then 
	$$|z|+\sum \limits _{n=2}^{\infty} (|a_{n}|+|b_{n}|) |z|^{n} \leq d(f(0),\partial f(\mathbb{D}))
	$$ for $|z|=r \leq r_{f}$, where $r_{f}$ is the unique  positive root of 
	\begin{equation} \label{himadri-vasu-p1-e-2.7}
	r+2 \sum \limits _{n=2}^{\infty}  \dfrac{r^{n}}{\alpha n^{2}+n(1-\alpha)}=1+2\sum \limits _{n=2}^{\infty} \dfrac{(-1)^{n-1}}{\alpha n^{2}+n(1-\alpha)} \end{equation} in $(0,1).$ 
	The radius $r_{f}$ is the best possible.
\end{thm}
\begin{rem}

\begin{enumerate}
\item [(i)] If $\alpha=0$,  the class $\mathcal{W}^{0}_{\mathcal{H}}(\alpha)$ reduces to $\mathcal{P}^{0}_{\mathcal{H}}.$ 
From \eqref{himadri-vasu-p1-e-2.7}, we obtain Bohr radius $r_{f} \approx0.285194$ for the class $\mathcal{P}^{0}_{\mathcal{H}}.$ 
\item [(ii)] If $\alpha=1,$ the class $\mathcal{W}^{0}_{\mathcal{H}}(\alpha)$ reduces to $\mathcal{W}^{0}_{\mathcal{H}}$ and we obtain the Bohr radius $r_{f}\approx 0.58387765$ for the class $\mathcal{W}^{0}_{\mathcal{H}}.$
\item[(iii)] When the co-analytic part $g\equiv 0,$ then $\mathcal{W}^{0}_{\mathcal{H}}(\alpha)$ reduces to $\mathcal{W}(\alpha).$ 
Therefore, from Lemma \ref{himadri-vasu-p1-lemma-002}, we  observe that Bohr radius for the class $\mathcal{W}(\alpha)$ is same as that of   the class $\mathcal{W}^{0}_{\mathcal{H}}(\alpha).$  
\end{enumerate}
\end{rem}
Using Lemmas \ref{himadri-vasu-p1-lemma-004} and \ref{himadri-vasu-p1-lemma-005}, we establish the Bohr phenomenon for the class $\mathcal{G}_{\mathcal{H}}^{k}(\alpha)$.
\begin{thm} \label{himadri-vasu-p1-theorem-004}
	Let $f=h+\overline{g} \in \mathcal{G}_{\mathcal{H}}^{k}(\alpha) $ be given by  \eqref{himadri-vasu-p1-e-1.7} and $\alpha \geq 0.$ Then $$
	|z|+\sum\limits_{n=k+1}^{\infty} (|a_{n}|+|b_{n}|)|z|^{n} \leq d(f(0),\partial f(\mathbb{D}))$$ for $|z|=r\leq r_{f}$, where $r_{f}$ is the smallest positive root of  
	\begin{equation} \label{himadri-vasu-p1-e-2.9}
	r+ 2\sum\limits_{n=k+1}^{\infty} \dfrac{r^{n}}{1+(n-1)\alpha}= 1+2\sum\limits_{n=1}^{\infty} \dfrac{(-1)^{n}}{1+nk\alpha} \end{equation} in $(0,1)$. The radius $r_{f}$ is the best possible.
\end{thm}
As an application of Lemma \ref{himadri-vasu-p1-lemma-008}, we obtain the sharp Bohr radius for the class $\mathcal{T}_{\mathcal{B}^{0}_{\mathcal{H}}(M)}.$
\begin{thm} \label{himadri-vasu-p1-theorem-005}
	Let $f \in \mathcal{T}_{\mathcal{B}^{0}_{\mathcal{H}}(M)}$ for $0< M <2$ be given by \eqref{himadri-vasu-p1-e-1.6}. Then 
	$$|z|+\sum\limits_{n=2}^{\infty} (|a_{n}|+|b_{n}|)|z|^{n}\leq d(f(0),\partial f(\mathbb{D}))$$
	for $|z|=r \leq r_{f}:=(-1+\sqrt{1+2M-M^{2}})/M,$ where $r_{f}$ is the positive root of 
	\begin{equation} \label{himadri-vasu-p1-e-2.10-a}
	M r^{2}+2r+(M-2)=0.\end{equation} The radius $r_{f}$ is the best possible.
\end{thm}

The Jacobian of a complex-valued harmonic function $f=h+\overline{g}$ is defined by 
$J_{f}=|h'(z)|^{2}- |g'(z)|^{2}$. It is known that if $J_{f}>0$ then the map $f$ is sense preserving and if $J_{f}<0$ then the map $f$ is sense reversing. 
Then using Lemma \ref{himadri-vasu-p1-lemma-008}, we prove the following theorem.
\begin{thm} \label{himadri-vasu-p1-theorem-006}
	Every function  $f \in \mathcal{T}_{\mathcal{B}^{0}_{\mathcal{H}}(M)}$ given by \eqref{himadri-vasu-p1-e-1.6} satisfies  
	$$|f(z)|+\sqrt{|J_{f}(z)|} \,\,|z| + \sum\limits_{n=2}^{\infty} (|a_{n}|+|b_{n}|)|z|^{n}\leq d(f(0),\partial f(\mathbb{D}))
	$$
	for $|z|=r\leq r_{f}:=(-1+\sqrt{1+2M-M^{2}})/2M$, where $r_{f}$ is the unique positive root of 
	$$
	4Mr^{2}+4r+(M-2)=0.
	$$ 
	in $(0,1)$. The radius $r_{f}$ is the best possible.
\end{thm}

 In the next theorem, we obtain the sharp left hand inequality in Theorem \ref{himadri-vasu-p1-theorem-006-b}, which will help us to obtain the Euclidean distance between $f(0)$ and $f(\mathbb{D}).$
\begin{thm} \label{himadri-vasu-p1-theorem-006-a}
Let $f \in \mathcal{P}^{0}_{\mathcal{H}}(M)$. Then 
\begin{equation} \label{himadri-vasu-p1-e-3.32}
|z| +2M \sum\limits_{n=2}^{\infty} \dfrac{(-1)^{n-1}|z|^{n}}{n(n-1)} \leq |f(z)| \leq |z| + 2M \sum\limits_{n=2}^{\infty} \dfrac{|z|^{n}}{n(n-1)}.
\end{equation}
Both the equalities are sharp for the function $f_{M}$ given by $f_{M}(z)=z+ 2M \sum\limits_{n=2}^{\infty} \dfrac{z^n}{n(n-1)} .
$
\end{thm}

Using Lemma \ref{himadri-vasu-p1-lemma-009} and Theorem \ref{himadri-vasu-p1-theorem-006-a}, we obtain the sharp Bohr radius for the class $\mathcal{P}^{0}_{\mathcal{H}}(M)$.
\begin{thm} \label{himadri-vasu-p1-theorem-007}
	Let $f \in \mathcal{P}^{0}_{\mathcal{H}}(M)$ be given by \eqref{himadri-vasu-p1-e-1.6} with $0<M< 1/(2(ln \, 4 -1))$. Then  
	$$|z|+\sum\limits_{n=2}^{\infty} (|a_{n}|+|b_{n}|)|z|^{n} \leq d(f(0),\partial f(\mathbb{D})) 
	$$
	for $|z|=r \leq r_{f} $,  where $r_{f}$ is the unique positive  root of 
	\begin{equation} \label{himadri-vasu-p1-e-2.13}
 r+ 2M \sum\limits_{n=2}^{\infty} \dfrac{r^n}{n(n-1)}=1+2M \sum\limits_{n=2}^{\infty} \dfrac{(-1)^{n-1}}{n(n-1)}
  \end{equation} 
  in $(0,1)$. The radius $r_{f}$ is the best possible.
\end{thm}
\begin{rem}
When the co-analytic part $g\equiv 0$ then  $\mathcal{P}^{0}_{\mathcal{H}}(M)$ coincides with $\mathcal{P}(M)$. 
Using Lemma \ref{himadri-vasu-p1-lemma-007}, it is easy to see that Bohr radius for the class $\mathcal{P}(M)$ is same as that of  the class $\mathcal{P}^{0}_{\mathcal{H}}(M).$
\end{rem}

\section{The Proof of main results}

\begin{pf}[{\bf Proof of Theorem   \ref{himadri-vasu-p1-theorem-001-a}}] 
	Let $f \in \mathcal{P}^{0}_{\mathcal{H}}(\alpha)$ then $$\real (h'(z)-\alpha)>|g'(z)| \quad \mbox{for} \quad z\in \mathbb{D},
	$$ which follows that $\real (h'(z)+\epsilon g'(z)-\alpha)>0$ for each $\epsilon \, (|\epsilon|=1)$ and $z \in \mathbb{D}$. In view of the subordination principle, there exists an analytic  function $\omega:\mathbb{D} \rightarrow \mathbb{D}$ with $\omega(0)=0$ such that $$\dfrac{h'(z)+\epsilon g'(z)-\alpha}{1-\alpha}=\dfrac{1+\omega(z)}{1-\omega(z)},$$ equivalently, 
	\begin{equation} \label{himadri-vasu-p1-e-3.1}
	h'(z)+\epsilon g'(z)=\alpha+(1-\alpha) \left(\dfrac{1+\omega(z)}{1-\omega(z)}\right).
	\end{equation}
	 Let  $F(z)=h(z)+\epsilon g(z)$ then from \eqref{himadri-vasu-p1-e-3.1}, it follows that 
	\begin{equation} \label{himadri-vasu-p1-e-3.2}
	F'(z)= \alpha+(1-\alpha) \left(\dfrac{1+\omega(z)}{1-\omega(z)} \right).
	\end{equation}
	Integrating \eqref{himadri-vasu-p1-e-3.2} along the linear segment connecting the origin and $z\in \mathbb{D}$, we obtain 
	\begin{equation} \label{himadri-vasu-p1-e-3.3}
	F(z)=\int \limits _{0}^{z} F'(\xi) \, d\xi  =\int \limits _{0}^{|z|} \left(\alpha+(1-\alpha) \left(\dfrac{1+\omega(te^{i\theta})}{1-\omega(te^{i\theta})}\right)\right) e^{i\theta} dt
	\end{equation}
	and hence \begin{align} \label{himadri-vasu-p1-e-3.3-a}
	|F(z)| & =\left|\int \limits _{0}^{|z|} \left(\alpha+(1-\alpha)\left(\dfrac{1+\omega(te^{i\theta})}{1-\omega(te^{i\theta})}\right)\right)dt\right| \\ \nonumber
	& \leq \int \limits _{0}^{|z|} \left(\alpha+(1-\alpha) \left(\dfrac{1+t}{1-t}\right)\right)dt \\ \nonumber
	&= |z|+(1-\alpha) \sum \limits _{n=2}^{\infty}  \dfrac{2}{n} |z|^{n}.
	\end{align}
	Similarly from \eqref{himadri-vasu-p1-e-3.3}, we obtain \begin{align} \label{himadri-vasu-p1-e-3.3-b}
	|F(z)|& \geq \int \limits _{0}^{|z|} \left(\alpha+(1-\alpha) \real \left(\dfrac{1+\omega(te^{i\theta})}{1-\omega(te^{i\theta})}\right)\right)dt \\ \nonumber & \geq 
	\int \limits _{0}^{|z|} \left(\alpha+(1-\alpha)\left(\dfrac{1-t}{1+t}\right)\right)dt \\ \nonumber & = 
	|z|+(1-\alpha) \sum \limits _{n=2}^{\infty}  \dfrac{2(-1)^{n-1}}{n} |z|^{n}.
	\end{align}
	From \eqref{himadri-vasu-p1-e-3.3-a} and \eqref{himadri-vasu-p1-e-3.3-b} we have $$|z|+ (1-\alpha) \sum\limits_{n=2}^{\infty}  \dfrac{2(-1)^{n-1}}{n} |z|^{n} \leq |F(z)| \leq |z|+ (1-\alpha)\sum\limits_{n=2}^{\infty}  \dfrac{2}{n} |z|^{n}.$$ Since $\epsilon \,(|\epsilon|=1)$ is arbitrary, for each $0\leq \alpha <1$, we have $$|z|+(1-\alpha) \sum\limits_{n=2}^{\infty}  \dfrac{2(-1)^{n-1}}{n} |z|^{n} \leq |f(z)| \leq |z|+(1-\alpha) \sum\limits_{n=2}^{\infty}  \dfrac{2}{n} |z|^{n}.$$ To establish the equality in \eqref{himadri-vasu-p1-e-2.2}, we consider the function $f_{\alpha}(z)$ defined by  $$f_{\alpha}(z)=z+\sum\limits_{n=2}^{\infty} \dfrac{2(1-\alpha)}{n} z^{n}. $$ It is easy to see that $f \in \mathcal{P}^{0}_{\mathcal{H}}(\alpha)$ for $0\leq \alpha < 1$. Let $|z|=r.$ The equality in both the sides of  \eqref{himadri-vasu-p1-e-2.2} holds for the function $f_{\alpha}$ at $z=-r$ and $z=r$ repectively.
\end{pf}

\begin{pf}[{\bf Proof of Theorem \ref{himadri-vasu-p1-theorem-001}}] 
Let $f \in \mathcal{P}^{0}_{\mathcal{H}}(\alpha)$ then from Theorem \ref{himadri-vasu-p1-theorem-001-a}, we have
	\begin{equation} \label{himadri-vasu-p1-e-3.4}
	|f(z)|\geq |z|+(1-\alpha)\sum\limits_{n=2}^{\infty}  \dfrac{2(-1)^{n-1}}{n} |z|^{n} \quad \mbox{for } \quad |z|<1.
	\end{equation}
	By taking $\liminf$ as $|z|\rightarrow 1$ on both sides of \eqref{himadri-vasu-p1-e-3.4}, we obtain 
	\begin{equation} \label{himadri-vasu-p1-e-3.5}
	\liminf \limits_{|z|\rightarrow 1} |f(z)| \geq 1+\sum\limits_{n=2}^{\infty}  2(1-\alpha) \dfrac{(-1)^{n-1}}{n}. 
	\end{equation}
	 The Euclidean distance between $f(0)$ and the boundary of $f(\mathbb{D})$ is given by \begin{equation} \label{himadri-vasu-p1-e-3.6}
	d(f(0), \partial f(\mathbb{D}))= \liminf \limits_{|z|\rightarrow 1} |f(z)-f(0)|.
	\end{equation}
	Since $f(0)=0$ from \eqref{himadri-vasu-p1-e-3.5} and \eqref{himadri-vasu-p1-e-3.6}, we obtain \begin{equation} \label{himadri-vasu-p1-e-3.7}
	d(f(0), \partial f(\mathbb{D})) \geq 1+\sum\limits_{n=2}^{\infty}  2(1-\alpha) \dfrac{(-1)^{n-1}}{n}.
	\end{equation}
	Let $H_{1}:[0,1) \rightarrow \mathbb{R}$ be defined by $$
	H_{1}(r)=r+2(1-\alpha)\sum \limits _{n=2}^{\infty}\dfrac{r^{n}}{n}-1-2(1-\alpha)\sum \limits _{n=2}^{\infty}\dfrac{(-1)^{n-1}}{n}.
	$$ 
	Clearly, $H_{1}$ is continuous in $[0,1)$ and differentiable in $(0,1).$
	Note that $$H_{1}(0)=-1-2(1-\alpha)\sum \limits _{n=2}^{\infty}\dfrac{(-1)^{n-1}}{n}$$
	and $\sum_{n=2}^{\infty}(-1)^{n-1}/n=ln\,2 -1.$ Hence $H_{1}(0)<0$ for each $\alpha \in [0,1).$ On the other hand, since $H_{1}(r)\rightarrow +\infty$ as $r \rightarrow 1$ and $$
	H'_{1}(r)=1+2(1-\alpha)\dfrac{r}{1-r}>0 $$
	for $r\in (0,1),$ it follows that $H_{1}(r)$ is strictly increasing in $(0,1).$\\
	Since $H_{1}(0)<0$ and $H_{1}(r)\rightarrow +\infty$ as $r \rightarrow 1$, the monotonocity of $H_{1}(r)$ implies that $H_{1}(r)$ has exactly one zero in $(0,1).$ 
	Let $r_{f}$ be the unique root of $H_{1}(r)$ in $(0,1).$ Then $H_{1}(r_{f})=0,$ which is equivalent to   
	\begin{equation} \label{himadri-vasu-p1-e-3.8}
	r_{f}+2(1-\alpha)\sum \limits_{n=2}^{\infty} \dfrac{r_{f}^{n}}{n}=1+\sum\limits_{n=2}^{\infty}  2(1-\alpha) \dfrac{(-1)^{n-1}}{n}.\end{equation}
	 For $0<r\leq r_{f}$, it follows from \eqref{himadri-vasu-p1-e-3.8} that \begin{equation} \label{himadri-vasu-p1-e-3.9}
	r+2(1-\alpha)\sum \limits _{n=2}^{\infty} \dfrac{r^{n}}{n}\leq r_{f}+2(1-\alpha)\sum \limits_{n=2}^{\infty} \dfrac{r_{f}^{n}}{n} = 1+\sum\limits_{n=2}^{\infty}  2(1-\alpha) \dfrac{(-1)^{n-1}}{n}.\end{equation} Using Lemma \ref{himadri-vasu-p1-lemma-001} and the inequalities \eqref{himadri-vasu-p1-e-3.7} and  \eqref{himadri-vasu-p1-e-3.9}, for $0<|z|=r\leq r_{f}$, we obtain \begin{align*}
	|z|+\sum\limits_{n=2}^{\infty}(|a_{n}|+|b_{n}|)|z|^{n}   & \leq r+2(1-\alpha)\sum \limits _{n=2}^{\infty}\dfrac{r^{n}}{n} \\ & \leq 1+\sum\limits_{n=2}^{\infty}  2(1-\alpha) \dfrac{(-1)^{n-1}}{n} \\ & \leq d(f(0), \partial f(\mathbb{D})).\end{align*}
	To show that the constant $r_{f}$ is the best possible, we consider the function $f_{\alpha}(z)$ defined by 
	$$f_{\alpha}(z)=z+\sum\limits_{n=2}^{\infty} \dfrac{2(1-\alpha)}{n} z^{n}.$$
	It is easy to see that the function $f_{\alpha} \in \mathcal{P}^{0}_{\mathcal{H}}(\alpha).$ 
	Let $|z|=r_{f}$, then a simple computation using \eqref{himadri-vasu-p1-e-3.8} gives 
	\begin{align*}
	|z|+\sum\limits_{n=2}^{\infty}(|a_{n}|+|b_{n}|)|z|^{n} &=r_{f}+2(1-\alpha)\sum \limits_{n=2}^{\infty} \dfrac{r_{f}^{n}}{n} \\ & = 1+\sum\limits_{n=2}^{\infty}  2(1-\alpha) \dfrac{(-1)^{n-1}}{n} \\ &=d(f(0), \partial f(\mathbb{D}))\end{align*} and hence the radius $r_{f}$ is the best possible. This completes the proof.
\end{pf}

\begin{pf}[{\bf Proof of Theorem  \ref{himadri-vasu-p1-theorem-002}}] 
	Let $f =h+\overline{g} \in \widetilde{\mathcal{G}}^{0}_{\mathcal{H}}(\beta)$ then  $$\real \left(\dfrac{h(z)}{z}\right)-\beta > \left|\dfrac{g(z)}{z} \right|,$$ which implies that 
	$$\real \left(\dfrac{h(z)}{z}+\epsilon \dfrac{g(z)}{z}-\beta \right) >0$$ for $z\in \mathbb{D}$ and each $\epsilon \, (|\epsilon|=1)$. In view of the subordination principle, there 
	exists an analytic function $p$ of the form $p(z)=1+\sum_{n=1}^{\infty} p_{n} z^{n}$ with $\real p(z)>0$ in $\mathbb{D}$ such that 
	\begin{equation} \label{himadri-vasu-p1-e-3.10}
	\dfrac{h(z)}{z}+\epsilon \dfrac{g(z)}{z}=\beta +(1-\beta)p(z).
	\end{equation}
	Comparing coefficients on both the sides of \eqref{himadri-vasu-p1-e-3.10}, we obtain 
	\begin{equation} \label{himadri-vasu-p1-e-3.11}
	a_{n}+\epsilon b_{n}=(1-\beta) p_{n-1} \quad \mbox{for} \quad n\geq 2. 
	\end{equation}
	Since $|p_{n}|\leq 2$ for $n\geq 1$ and $\epsilon \, (|\epsilon|=1)$ is arbitrary, it follows from \eqref{himadri-vasu-p1-e-3.11} that 
	\begin{equation} \label{himadri-vasu-p1-e-3.12}
	|a_{n}|+|b_{n}| \leq 2(1-\beta).
	\end{equation}
	Let $F(z)=h(z)+\epsilon g(z)$. Then from \eqref{himadri-vasu-p1-e-3.10}, we obtain \begin{align*}
	\left|\dfrac{F(z)}{z} \right|&= |\beta +(1-\beta)p(z)| \\ &\leq \beta +(1-\beta) \left(\dfrac{1+|z|}{1-|z|}\right),
	\end{align*}
	which is equivalent to \begin{equation} \label{himadri-vasu-p1-e-3.13}
	|F(z)|\leq \beta |z|+(1-\beta) \left(\dfrac{1+|z|}{1-|z|}\right) |z|. 
	\end{equation}
	Again, from \eqref{himadri-vasu-p1-e-3.10} we obtain \begin{align*}
	\left|\dfrac{F(z)}{z}\right| &\geq \beta +(1-\beta)\real p(z) \\ & \geq \beta +(1-\beta) \left(\dfrac{1-|z|}{1+|z|}\right),
	\end{align*}
	which is equivalent to \begin{equation} \label{himadri-vasu-p1-e-3.14}
	|F(z)|\geq \beta |z|+(1-\beta)\left(\dfrac{1-|z|}{1+|z|}\right) |z|.
	\end{equation}
	Therefore from \eqref{himadri-vasu-p1-e-3.13} and \eqref{himadri-vasu-p1-e-3.14}, we obtain
	\begin{equation} \label{himadri-vasu-p1-e-3.15}
	\beta |z|+(1-\beta) \left(\dfrac{1-|z|}{1+|z|}\right) |z| \leq |F(z)| \leq \beta |z| +(1-\beta) \left(\dfrac{1+|z|}{1-|z|}\right) |z|.\end{equation}
	Since $F(z)=h(z)+\epsilon g(z)$ and $\epsilon(|\epsilon|=1)$ is arbitrary, \eqref{himadri-vasu-p1-e-3.15} yields  \begin{equation}  \label{himadri-vasu-p1-e-3.16}
	\beta |z|+(1-\beta)\left(\dfrac{1-|z|}{1+|z|}\right) |z| \leq |f(z)| \leq \beta |z| +(1-\beta) \left(\dfrac{1+|z|}{1-|z|}\right)|z|.
	\end{equation} 
	The equality holds for both the sides of \eqref{himadri-vasu-p1-e-3.16} for the following function \begin{equation} \label{himadri-vasu-p1-e-3.17}
	f(z)=z+\sum\limits_{n=2}^{\infty} 2(1-\beta) z^{n}
	\end{equation} and its rotations.
	From the left side inequality of \eqref{himadri-vasu-p1-e-3.16}, we have
	\begin{equation} \label{himadri-vasu-p1-e-3.18}
	|f(z)|\geq \beta |z|+(1-\beta) \left(\dfrac{1-|z|}{1+|z|}\right) |z|.\end{equation} 
	Taking $\liminf$ as $ |z|\rightarrow 1$ on both the sides of \eqref{himadri-vasu-p1-e-3.18}, we obtain \begin{equation}
	\liminf \limits_{|z|\rightarrow 1} |f(z)|\geq \beta.
	\end{equation} Therefore, the Euclidean distance between $f(0)$ and the boundary of $f(\mathbb{D})$ is 
	\begin{equation} \label{himadri-vasu-p1-e-3.20}
	d(f(0),\partial f(\mathbb{D}))=\liminf \limits_{|z|\rightarrow 1} |f(z)-f(0)| \geq \beta.
	\end{equation}
	For $|z|=r$, using \eqref{himadri-vasu-p1-e-3.12} we obtain 
	\begin{equation} \label{himadri-vasu-p1-e-3.20-b}
	|z|+\sum\limits_{n=2}^{\infty}(|a_{n}|+|b_{n}|)|z|^{n} \leq r+2(1-\beta) \sum\limits_{n=2}^{\infty} r^{n}=	r+2(1-\beta) \dfrac{r^{2}}{1-r}.\end{equation}	
	A simple computation shows that 
	\begin{equation} \label{himadri-vasu-p1-e-3.20-c}
	r+2(1-\beta) \dfrac{r^{2}}{1-r}\leq \beta 
	\end{equation}
	when $$
	H(r):=r+2(1-\beta) \dfrac{r^{2}}{1-r}-\beta \leq 0.$$ 
     Note that $H(r)=0$ is equivalent to 
     \begin{equation} \label{himadri-vasu-p1-e-3.20-a}
     (1-2\beta)r^{2}+(1+\beta)r-\beta =0 
     \end{equation}
     and has the following roots  
	$$
	r=\dfrac{-1-\beta \pm \sqrt{1+6\beta -7\beta ^{2}}}{2(1-2\beta)}.$$
	We note  that $r=(-1-\beta +\sqrt{1+6\beta -7\beta ^{2}})/2(1-2\beta)$ is the only positive root of (\ref{himadri-vasu-p1-e-3.20-a})  which lies in $(0,1)$ for $0\leq \beta <1/2$ and we choose this to be $r_{f}.$ 
	So $r_{f}$ satisfies the equation \eqref{himadri-vasu-p1-e-3.20-a}, which in turns equivalent to  \begin{equation} \label{himadri-vasu-p1-e-3.21}
	r_{f}+2(1-\beta) \dfrac{r^{2}_{f}}{1-r_{f}}=\beta.
	\end{equation} 
	Clearly, $H(r)$ is continuous in $[0,1)$ and differentiable in $(0,1)$. A simple calculation shows that $$
	H'(r)=1+2\, (1-\beta) \dfrac{r(2-r)}{(1-r)^{2}}>0$$
	for $r\in(0,1)$ and $0\leq \beta <1/2$. Thus, $H(r)$ is strictly increasing in $(0,1)$ for $0\leq \beta <1/2$. Therefore, for $r\leq r_{f}$, we have $H(r) \leq H(r_{f})=0$. 
	So we conclude that the inequality \eqref{himadri-vasu-p1-e-3.20-c} is satisfied if $r\leq r_{f}$.\\
	For $0< |z|=r \leq r_{f}$, it follows from  \eqref{himadri-vasu-p1-e-3.20}, \eqref{himadri-vasu-p1-e-3.20-b}, \eqref{himadri-vasu-p1-e-3.20-c} and \eqref{himadri-vasu-p1-e-3.21} that 
	\begin{align*}
	|z|+\sum\limits_{n=2}^{\infty}(|a_{n}|+|b_{n}|)|z|^{n} &\leq r+2(1-\beta) \dfrac{r^{2}}{1-r} \\ &\leq
	r_{f}+2(1-\beta) \dfrac{r^{2}_{f}}{1-r_{f}}\\ &=\beta \\ &\leq d(f(0),\partial f(\mathbb{D})).
	\end{align*}
	Let $f$ be given by \eqref{himadri-vasu-p1-e-3.17}, then for $|z|=r_{f}$, a simple calculation using \eqref{himadri-vasu-p1-e-3.21} shows that
	 \begin{align*} 
	|z|+\sum\limits_{n=2}^{\infty}(|a_{n}|+|b_{n}|)|z|^{n} &=r_{f}+2(1-\beta) \dfrac{r^{2}_{f}}{1-r_{f}} \\ &=\beta \\ &=d(f(0),\partial f(\mathbb{D})) .\end{align*}
	This shows that the radius $r_{f}$ is the best possible. 
\end{pf}

\begin{pf}[{\bf Proof of Theorem \ref{himadri-vasu-p1-theorem-003}}] 
	Let $f \in \mathcal{W}^{0}_{\mathcal{H}}(\alpha)$  be given by \eqref{himadri-vasu-p1-e-1.6}. Then from Lemma \ref{himadri-vasu-p1-lemma-003-a}, it is evident that the Euclidean distance between $f(0)$ and the boundary of $f(\mathbb{D})$ is
	\begin{equation} \label{himadri-vasu-p1-e-3.22}
	d(f(0),\partial f(\mathbb{D})) =\liminf \limits_{|z|\rightarrow 1} |f(z)-f(0)| \geq 1+2\sum\limits_{n=2}^{\infty} \dfrac{(-1)^{n-1}}{\alpha n^{2}+n(1-\alpha)}.
	\end{equation}
	Let $H_{2}:[0,1] \rightarrow \mathbb{R}$ be defined by $$H_{2}(r)=r+ \sum \limits _{n=2}^{\infty}  \dfrac{2r^{n}}{\alpha n^{2}+n(1-\alpha)}- 1-2\sum \limits _{n=2}^{\infty} 
		\dfrac{(-1)^{n-1}}{\alpha n^{2}+n(1-\alpha)}.$$
	Clearly, $H_{2}$ is continuous on $[0,1]$ and differentiable on $(0,1)$ and
		$$H_{2}(0)=-1-2\sum \limits _{n=2}^{\infty} 
		\dfrac{(-1)^{n-1}}{\alpha n^{2}+n(1-\alpha)}.$$ Since $$\left|\sum \limits _{n=2}^{\infty} 
		\dfrac{(-1)^{n-1}}{\alpha n^{2}+n(1-\alpha)}\right| \leq \dfrac{1}{2} \quad \mbox{for} \quad n\geq 2,$$ we have $H_{2}(0)<0.$ On the other hand, 
		$$H_{2}(1)=2\sum \limits _{n=2}^{\infty}  \dfrac{1}{\alpha n^{2}+n(1-\alpha)} - 2\sum \limits _{n=2}^{\infty}\dfrac{(-1)^{n-1}}{\alpha n^{2}+n(1-\alpha)}$$ and $$2\sum \limits _{n=2}^{\infty}  \dfrac{1}{\alpha n^{2}+n(1-\alpha)} > 2\sum \limits _{n=2}^{\infty} \dfrac{(-1)^{n-1}}{\alpha n^{2}+n(1-\alpha)} \quad \mbox{for} \quad n\geq 2.$$ This shows that $H_{2}(1)>0.$ Since $H_{2}(0)<0$ and $H_{2}(1)>0$, by the intermediate value theorem, we conclude that $H_{2}(r)$ has a real root in $(0,1).$ To show that $H(r)$ has exactly one zero in $(0,1),$ it is enough to show that $H_{2}$ is monotonic on $(0,1).$ Since $$
		H'_{2}(r)=1+\sum \limits _{n=2}^{\infty}  \dfrac{2 n r^{n-1}}{\alpha n^{2}+n(1-\alpha)} > 0$$ for $r\in (0,1),$ $H_{2}(r)$ is strictly monotonically increasing in $(0,1).$ Therefore, $H_{2}(r)$ has exactly one zero in $(0,1).$ Let $r_{f}$ be the unique root of $H_{2}(r)$ in $(0,1).$
	Then $r_{f}$ satisfies 
	\begin{equation} \label{himadri-vasu-p1-e-3.23}
	r_{f}+ \sum \limits _{n=2}^{\infty}  \dfrac{2r_{f}^{n}}{\alpha n^{2}+n(1-\alpha)}=1+2\sum \limits _{n=2}^{\infty} 
	\dfrac{(-1)^{n-1}}{\alpha n^{2}+n(1-\alpha)}.
	\end{equation}  
	Using Lemma \ref{himadri-vasu-p1-lemma-003}, \eqref{himadri-vasu-p1-e-3.22} and \eqref{himadri-vasu-p1-e-3.23} for $0<|z|=r\leq r_{f}$, we obtain 
	\begin{align*}
	|z|+\sum \limits _{n=2}^{\infty} (|a_{n}|+|b_{n}|) |z|^{n} &\leq r+ \sum \limits _{n=2}^{\infty} \dfrac{2}{\alpha n^{2}+n(1-\alpha)} r^{n}\\ 
	&\leq r_{f}+ \sum \limits _{n=2}^{\infty} \dfrac{2}{\alpha n^{2}+n(1-\alpha)} r_{f}^{n} \\ &\leq d(f(0),\partial f(\mathbb{D})).
	\end{align*}
	To show that $r_{f}$ is the best possible, we consider the following function $f_{\alpha}$ defined by $$
	f_{\alpha}(z)=z+\sum \limits _{n=2}^{\infty} \dfrac{2}{\alpha n^{2}+n(1-\alpha)} z^{n}.$$ 
	Clearly, $f_{\alpha}$ belongs to the class $\mathcal{W}^{0}_{\mathcal{H}}(\alpha).$ For $f=f_{\alpha}$, we have \begin{equation} \label{himadri-vasu-p1-e-3.23-a}
	d(f(0),\partial f(\mathbb{D}))=1+2\sum \limits _{n=2}^{\infty} 
	\dfrac{(-1)^{n-1}}{\alpha n^{2}+n(1-\alpha)}.
	\end{equation}
	Using \eqref{himadri-vasu-p1-e-3.23} and \eqref{himadri-vasu-p1-e-3.23-a}, a simple computation with $|z|=r_{f}$ shows that \begin{align*}
	|z|+\sum \limits _{n=2}^{\infty} (|a_{n}|+|b_{n}|) |z|^{n} &= r_{f}+ \sum \limits _{n=2}^{\infty} \dfrac{2}{\alpha n^{2}+n(1-\alpha)} r_{f}^{n} \\ &= 1+2\sum \limits _{n=2}^{\infty} 
	\dfrac{(-1)^{n-1}}{\alpha n^{2}+n(1-\alpha)} \\ &=  d(f(0),\partial f(\mathbb{D})).
	\end{align*} Therefore the radius $r_{f}$ is the best possible.
\end{pf}

\begin{pf}[{\bf Proof of Theorem \ref{himadri-vasu-p1-theorem-004}}] 
	From the left side inequality in Lemma \ref{himadri-vasu-p1-lemma-005}, it is evident that the Euclidean distance between $f(0)$ and the boundary of $f(\mathbb{D})$ is
	\begin{equation} \label{himadri-vasu-p1-e-3.24}
	d(f(0),\partial f(\mathbb{D}))=\liminf \limits_{|z|\rightarrow 1} |f(z)-f(0)| \geq 1+2\sum\limits_{n=1}^{\infty} \dfrac{(-1)^{n}}{1+nk \alpha}.
	\end{equation}
	Let $H_{3}:[0,1)\rightarrow \mathbb{R}$ be defined by 
	$$	r+ 2\sum \limits _{n=k+1}^{\infty}  \dfrac{r^{n}}{1+(n-1)\alpha}-1-2\sum \limits _{n=1}^{\infty} \dfrac{(-1)^{n}}{1+nk \alpha}.
	$$ Then going by the same lines of argument as in the proof of 
	Theorem \ref{himadri-vasu-p1-theorem-001} and Theorem \ref{himadri-vasu-p1-theorem-003}, we obtain $H_{3}(r)$ has exactly one zero in $(0,1)$ and we choose it to be $r_{f}.$
	Since $r_{f}$ is a root of $H_{3}(r)$, we have $H_{3}(r_{f})=0$. Therefore  
	\begin{equation} \label{himadri-vasu-p1-e-3.25}
	r_{f}+ 2\sum \limits _{n=k+1}^{\infty}  \dfrac{r_{f}^{n}}{1+(n-1)\alpha}=1+2\sum \limits _{n=1}^{\infty} \dfrac{(-1)^{n}}{1+nk \alpha}. 
	\end{equation} 
	In view of Lemma \ref{himadri-vasu-p1-lemma-004}, \eqref{himadri-vasu-p1-e-3.24} and \eqref{himadri-vasu-p1-e-3.25}, for $0<|z|=r\leq r_{f}$, we obtain 
	\begin{align*}
	|z|+\sum \limits _{n=k+1}^{\infty} (|a_{n}|+|b_{n}|) |z|^{n} &\leq r+ 2\sum\limits_{n=k+1}^{\infty} \dfrac{r^{n}}{1+(n-1)\alpha}\\ &\leq r_{f}
	+ 2\sum \limits _{n=k+1}^{\infty} \dfrac{r_{f}^{n}}{1+(n-1)\alpha} \\ & = 1+2\sum \limits _{n=1}^{\infty} \dfrac{(-1)^{n}}{1+nk \alpha} \\ &\leq d(f(0),\partial f(\mathbb{D})).
	\end{align*}
	In order to show that $r_{f}$ is the best possible constant, we consider the following function $f_{\alpha}$ by 
	$$
	f_{\alpha}(z)=z+ \sum\limits_{n=1}^{\infty} \dfrac{2}{1+nk \alpha } z^{nk+1}.
	$$ 
	It is not difficult to see that $f_{\alpha} \in \mathcal{G}_{\mathcal{H}}^{k}(\alpha)$. For $f=f_{\alpha}$, we have \begin{equation} \label{himadri-vasu-p1-e-3.25-a}
	d(f(0),\partial f(\mathbb{D}))=1+2\sum \limits _{n=1}^{\infty} \dfrac{(-1)^{n}}{1+nk \alpha}.
	\end{equation} 
	For $f=f_{\alpha}$ and $|z|=r_{f},$ a simple computation using \eqref{himadri-vasu-p1-e-3.25} and \eqref{himadri-vasu-p1-e-3.25-a} shows that \begin{align*}
	|z|+\sum\limits_{n=k+1}^{\infty} (|a_{n}|+|b_{n}|)|z|^{n} &= r_{f}+ 2\sum \limits _{n=k+1}^{\infty} 
	\dfrac{r_{f}^{n}}{1+(n-1)\alpha} \\ &= 1+2\sum \limits _{n=1}^{\infty} \dfrac{(-1)^{n}}{1+nk \alpha} \\ &=d(f(0),\partial f(\mathbb{D})).
	\end{align*} 
	Therefore the radius $r_{f}$ is the best possible. This completes the proof.
\end{pf}

\begin{pf}[{\bf Proof of Theorem \ref{himadri-vasu-p1-theorem-005}}] 
	From the left side inequality in Lemma \ref{himadri-vasu-p1-lemma-008}, it follows that the Euclidean distance between $f(0)$ and the boundary of $f(\mathbb{D})$ is given by
	\begin{equation} \label{himadri-vasu-p1-e-3.26}
	d(f(0), \partial f(\mathbb{D}))=\liminf \limits_{|z|\rightarrow 1} |f(z)-f(0)| \geq 1-\dfrac{M}{2}.
	\end{equation}
	Observe that roots of the quadratic equation \eqref{himadri-vasu-p1-e-2.10-a} are 
	$$	r=\dfrac{-1\pm\sqrt{1+2M-M^{2}}}{M}
	$$ 
	and  $r=(-1+\sqrt{1+2M-M^{2}})/M$ is the only root which  lies in $(0,1)$ for $0<M<2.$ Let $r_{f}$ be that root.
	Since $r_{f}=(-1+\sqrt{1+2M-M^{2}})/M$ is a root of \eqref{himadri-vasu-p1-e-2.10-a}, it satisfies  
	\begin{equation} \label{himadri-vasu-p1-e-3.27}
	r_{f} + \dfrac{M}{2} r_{f}^{2}=1-\dfrac{M}{2}.
	\end{equation}
	For $0<r \leq r_{f}$, it follows from \eqref{himadri-vasu-p1-e-3.27} that
	\begin{equation} \label{himadri-vasu-p1-e-3.28}
	r + \dfrac{M}{2} r^{2}\leq r_{f} + \dfrac{M}{2} r_{f}^{2}=1-\dfrac{M}{2}.
	\end{equation}
	
	We note that  $\{n(n-1)\}$ is an increasing sequence of positive real numbers for $n\geq 2$ and $n(n-1) \geq 2$ for $n \geq 2$
	and hence 
	\begin{equation} \label{himadri-vasu-p1-e-3.29}
	\sum\limits_{n=2}^{\infty} (|a_{n}|+|b_{n}|) \leq \dfrac{M}{2}.
	\end{equation}
	Using \eqref{himadri-vasu-p1-e-3.26}, \eqref{himadri-vasu-p1-e-3.28} and \eqref{himadri-vasu-p1-e-3.29}, for $|z|=r\leq r_{f},$ we obtain 
	$$
	|z|+\sum\limits_{n=2}^{\infty} (|a_{n}|+|b_{n}|)|z|^{n}\leq r+\dfrac{M}{2}r^{2} \leq 1-\dfrac{M}{2} = d(f(0),\partial f(\mathbb{D})).
	$$
	To show that the radius $r_{f}$ is the best possible, we define the function 
	$$ f(z)=z+\dfrac{M}{2} z^{2}.
	$$
	It is easy to see that $f \in \mathcal{T}_{\mathcal{B}^{0}_{\mathcal{H}}(M)}$. For  $|z|=r_{f}$, we have
	\begin{align*}|z|+\sum\limits_{n=2}^{\infty} (|a_{n}|+|b_{n}|)|z|^{n} &=r_{f}+\dfrac{M}{2}r_{f}^{2}\\ &=1-\dfrac{M}{2} \\ &=d(f(0), \partial f(\mathbb{D}))
	\end{align*} and hence the constant $r_{f}$ is the best possible.
	This completes the proof.	
\end{pf}

\begin{pf}[{\bf Proof of Theorem \ref{himadri-vasu-p1-theorem-006}}] 
	The Jacobian of a complex-valued harmonic function $f=h+\overline{g}$ is defined by 
	$J_{f}=|h'(z)|^{2}- |g'(z)|^{2}.$
	 If $f \in \mathcal{T}_{\mathcal{B}^{0}_{\mathcal{H}}(M)}$ then $\mathcal{T}_{\mathcal{B}^{0}_{\mathcal{H}}(M)} \subseteq \mathcal{B}^{0}_{\mathcal{H}}(M)$ 
	 and $J_{f}(z) \leq (1+M|z|)^{2}$ (see \cite{nirupam cvee}). By similar argument as in Theorem \ref{himadri-vasu-p1-theorem-005}, we can show the radius $r_{f}=(-1+\sqrt{1+2M-M^{2}})/2M$ in 
	 Theorem \ref{himadri-vasu-p1-theorem-006} satisfies 
	 \begin{equation} \label{himadri-vasu-p1-e-3.30}
	2Mr_{f}^{2}+2r_{f}=1-\dfrac{M}{2}.
	\end{equation}
	For $0<r \leq r_{f}$, it follows from \eqref{himadri-vasu-p1-e-3.30} that 
	\begin{equation} \label{himadri-vasu-p1-e-3.31}
	2Mr^{2}+2r \leq 2Mr_{f}^{2}+2r_{f}=1-\dfrac{M}{2}.
	\end{equation}
	A simple computation using Lemma \ref{himadri-vasu-p1-lemma-008}, \eqref{himadri-vasu-p1-e-3.29} and \eqref{himadri-vasu-p1-e-3.31} shows that 
	\begin{align*}
	|f(z)|+\sqrt{|J_{f}(z)|} \,\,|z| + \sum\limits_{n=2}^{\infty} (|a_{n}|+|b_{n}|)|z|^{n} 
	& \leq r+\dfrac{M}{2}r^{2} + (1+Mr)r + \dfrac{M}{2}r^{2} \\ 
	& \leq  2Mr^{2}+2r \\ &\leq  1-\dfrac{M}{2} \\  
	& \leq  d(f(0),\partial f(\mathbb{D})).
	\end{align*}
	To show the sharpness of the radius $r_{f}$, we consider the following function $$f(z)=z+\dfrac{M}{2} z^{2}.$$
	For $|z|=r_{f}$, a simple computation shows that 
	\begin{align*}|f(z)|+\sqrt{|J_{f}(z)|} \,\,|z| + \sum\limits_{n=2}^{\infty} (|a_{n}|+|b_{n}|)|z|^{n} &= 2Mr_{f}^{2}+2r_{f}\\ &= 1-\dfrac{M}{2} \\ &=d(f(0),\partial f(\mathbb{D}))
	\end{align*} 
	and hence the radius  $r_{f}$ is the best possible.
	This completes the proof.
\end{pf}

\begin{pf} [{\bf Proof of Theorem \ref{himadri-vasu-p1-theorem-006-a}}] 
Since sharp right hand side inequality of \eqref{himadri-vasu-p1-e-3.32} has already been established in \cite{nirupam-bull aus math}, we omit the proof. We prove the sharp left hand inequality.
Let $f =h+\overline{g}\in \mathcal{P}^{0}_{\mathcal{H}}(M)$. Then in view of Lemma \ref{himadri-vasu-p1-lemma-007}, $F_{\epsilon}=h+ \epsilon g$ belongs to $\mathcal{P}(M)$ for each $\epsilon$ with $|\epsilon|=1$ and hence we have $$
\real(zF''_{\epsilon}(z))=z(h''(z)+\epsilon g''(z))>-M \quad \mbox{for} \quad z \in \mathbb{D}.$$ 	
Therefore by the subordination principle, there exists an analytic function $\omega:\mathbb{D}\rightarrow \mathbb{D}$ with $\omega(0)=0$ such that \begin{equation} \label{himadri-vasu-p1-e-3.31-a}
\dfrac{zF''_{\epsilon}(z)+M}{M}=\dfrac{1+\omega(z)}{1-\omega(z)}
\end{equation} which is equivalent to 
\begin{equation} \label{himadri-vasu-p1-e-3.31-b}
F''_{\epsilon}(z)=\dfrac{2M\omega(z)}{z(1-\omega(z))}.
\end{equation}
Note that $|\omega(z)|\leq |z|$ and 
\begin{align} \label{himadri-vasu-p1-e-3.34}
\real \left(\dfrac{\omega(z)}{1-\omega(z)}\right) & =
\dfrac{1}{2}\left(\dfrac{\omega(z)}{1-\omega(z)} + \dfrac{\overline{\omega(z)}}{1-\overline{\omega(z)}}\right) \\ \nonumber & = \dfrac{1}{2}\left(\dfrac{2\real\omega(z)-2|\omega(z)|^{2}}{|1-\omega(z)|^{2}}\right) \\ \nonumber & \geq \dfrac{-|\omega(z)|}{1+|\omega(z)|} \\ \nonumber & \geq \dfrac{-|z|}{1+|z|}.
\end{align}
By integrating \eqref{himadri-vasu-p1-e-3.31-b} along line segment joining $0$ to $z$ and using  \eqref{himadri-vasu-p1-e-3.34}, we obtain 
\begin{align*}
|F'_{\epsilon}(z)|&= \left|1+2M \int_{0}^{z} \dfrac{\omega(\xi)}{\xi (1-\omega(\xi))}\,\, d\xi\right| \\ & \geq 1+ 2M \int \limits _{0}^{|z|} \dfrac{1}{t} \real \left(\dfrac{\omega (t e^{i\theta})}{ (1-\omega (t e^{i\theta}))}\right) \,dt \\ & \geq 1+ 2M \int \limits _{0}^{|z|} \dfrac{1}{t} \dfrac{-t}{1+t} \,dt \\ & = 1+ 2M \sum\limits_{n=2}^{\infty} (-1)^{n-1} \dfrac{|z|^{n-1}}{n-1}  
\end{align*}
	and hence \begin{equation} \label{himadri-vasu-p1-e-3.35}
|F'_{\epsilon}(z)| \geq 1+ 2M \sum\limits_{n=2}^{\infty} (-1)^{n-1} \dfrac{|z|^{n-1}}{n-1}.
\end{equation} Since $\epsilon \,(|\epsilon|=1)$ is arbitrary, it follows from \eqref{himadri-vasu-p1-e-3.35} that \begin{equation} \label{himadri-vasu-p1-e-3.36}
|h'(z)|-|g'(z)| \geq 1+ 2M \sum\limits_{n=2}^{\infty} (-1)^{n-1} \dfrac{|z|^{n-1}}{n-1}.
\end{equation} In view of \eqref{himadri-vasu-p1-e-3.36}, we obtain \begin{align} \label{himadri-vasu-p1-e-3.37}
|f(z)| &= \left|\int \limits _{0}^{z}  \dfrac{\partial f}{\partial \xi} d\xi + \dfrac{\partial f}{\partial \overline{\xi}} d \overline{\xi}\right| \\ \nonumber & \geq \int \limits _{0}^{z} (|h'(\xi)|-|g'(\xi)|) |d\xi| \\ \nonumber & \geq \int \limits _{0}^{|z|} \left(1+ 2M \sum\limits_{n=2}^{\infty} (-1)^{n-1} \dfrac{t^{n-1}}{n-1}\right) dt \\ \nonumber &= |z|+ 2M \sum\limits_{n=2}^{\infty} (-1)^{n-1} \dfrac{|z|^{n}}{n(n-1)}.
\end{align}
The equality in both sides of \eqref{himadri-vasu-p1-e-3.32} holds for the function $f=f_{M}$ given by $$
f_{M}(z)=z+ 2M \sum\limits_{n=2}^{\infty} \dfrac{z^n}{n(n-1)}$$ at $z=-r$ and $z=r$ respectively.
\end{pf}
\begin{pf}[{\bf Proof of Theorem \ref{himadri-vasu-p1-theorem-007}}] 
	Let $f=h+\overline{g} \in \mathcal{P}^{0}_{\mathcal{H}}(M)$ be given by \eqref{himadri-vasu-p1-e-1.6} for $0<M<1/(2(ln \, 4 -1)).$ 
		From \eqref{himadri-vasu-p1-e-3.37}, it is evident that the Euclidean distance between $f(0)$ and $f(\mathbb{D})$ is \begin{equation} \label{himadri-vasu-p1-e-3.38}
	d(f(0),\partial f(\mathbb{D})) \geq 1+2M \sum\limits_{n=2}^{\infty}  \dfrac{(-1)^{n-1}}{n(n-1)}.
	\end{equation}
	Let $H_{4}:[0,1] \rightarrow \mathbb{R}$ defined by $$H_{4}(r)=  r+ 2M \sum\limits_{n=2}^{\infty} \dfrac{r^n}{n(n-1)}-1- 2M \sum\limits_{n=2}^{\infty} \dfrac{(-1)^{n-1}}{n(n-1)}.$$ Clearly, $H_{4}$ is continuous in $[0,1]$ and differentiable in $(0,1)$.
	We note that $$H_{4}(0)=-1-2M \sum\limits_{n=2}^{\infty} \dfrac{(-1)^{n-1}}{n(n-1)}= -1-2M(1-ln \, 4) <0$$ when $0<M<1/(2(ln \, 4 -1))$ and  \begin{align*}
	H_{4}(1)&=2M \sum\limits_{n=2}^{\infty} \dfrac{1}{n(n-1)}- 2M \sum\limits_{n=2}^{\infty} \dfrac{(-1)^{n-1}}{n(n-1)}\\ &= 2M-  2M \sum\limits_{n=2}^{\infty} \dfrac{(-1)^{n-1}}{n(n-1)} \\ & = 2M(2- ln\,4)>0.\end{align*}
	Since for $0<M<1/(2(ln \, 4 -1))$, $H_{4}(0)<0$ and $H_{4}(1)>0$, by intermediate value theorem, we conclude that $H$ has a root in $(0,1).$ In order to prove the uniqueness of that root, it is enough to show that $H$ is monotonic function in $\mathbb{D}.$ Since $$H_{4}'(r)=1+2M \sum\limits_{n=2}^{\infty} \dfrac{r^{n-1}}{n-1}>0$$ for each $r\in (0,1)$, $H_{4}$ is strictly monotonic increasing in $\mathbb{D}$. Hence $H_{4}$ has exactly one zero in $(0,1)$ and let that root be $r_{f}$. Since $r_{f}$  is a root of \eqref{himadri-vasu-p1-e-2.13}, we have 
	\begin{equation} \label{himadri-vasu-p1-e-3.39}
	r_{f}+ 2M \sum\limits_{n=2}^{\infty} \dfrac{r_{f}^n}{n(n-1)}=1+ 2M \sum\limits_{n=2}^{\infty} \dfrac{(-1)^{n-1}}{n(n-1)}.  
	\end{equation}
	Using \eqref{himadri-vasu-p1-e-3.38}, \eqref{himadri-vasu-p1-e-3.39} and Lemma \ref{himadri-vasu-p1-lemma-009}, we obtain   
	\begin{align*}
	|z|+\sum\limits_{n=2}^{\infty} (|a_{n}|+|b_{n}|)|z|^{n} 
	&\leq  r+ 2M \sum\limits_{n=2}^{\infty} \dfrac{r^n}{n(n-1)}\\ 
	& \leq r_{f}+ 2M \sum\limits_{n=2}^{\infty} \dfrac{r_{f}^n}{n(n-1)} \\ &= 1+ 2M \sum\limits_{n=2}^{\infty} \dfrac{(-1)^{n-1}}{n(n-1)} 
	\\ & \leq d(f(0),\partial f(\mathbb{D}))
	\end{align*} 
	for $0<|z|=r \leq r_{f}$.
	It is not difficult to show that the function $f$ defined by 
	$$f_{M}(z)=z+ 2M \sum\limits_{n=2}^{\infty} \dfrac{z^n}{n(n-1)} 
	$$ 
	belongs to $\mathcal{P}^{0}_{\mathcal{H}}(M)$. 
	For $f=f_{M}$ and $|z|=r \leq r_{f}$, a simple computation shows that 
	\begin{align*}	|z|+\sum\limits_{n=2}^{\infty} (|a_{n}|+|b_{n}|)|z|^{n} &=r_{f}+2M\sum\limits_{n=2}^{\infty} \dfrac{r_{f}^n}{n(n-1)} \\ &= 1+ 2M \sum\limits_{n=2}^{\infty} \dfrac{(-1)^{n-1}}{n(n-1)} \\ &=d(f(0),\partial f(\mathbb{D}))
	\end{align*}
	and hence the radius $r_{f}$ is the best possible.
\end{pf}

\vspace{1cm}
\noindent\textbf{Acknowledgement:}  The first author thank SERB-MATRICS and the second author thank CSIR for their support. 

\end{document}